\documentclass[a4paper,11pt]{article}
\usepackage{amsmath, amsfonts, amscd, amssymb, amsthm, enumerate, xypic}

\newcommand{\lrk}{\operatorname{lrk}}
\newcommand{\Mat}{\operatorname{M}}
\newcommand{\Mata}{\operatorname{A}}

\newcommand{\GL}{\operatorname{GL}}

\newcommand{\Vect}{\operatorname{span}}

\newcommand{\tr}{\operatorname{tr}}

\newcommand{\rk}{\operatorname{rk}}
\newcommand{\codim}{\operatorname{codim}}
\renewcommand{\setminus}{\smallsetminus}

\def\K{\mathbb{K}}

\def\calE{\mathcal{E}}

\def\calH{\mathcal{H}}

\def\calV{\mathcal{V}}
\def\calW{\mathcal{W}}
\def\calX{\mathcal{X}}

\def\lcro{\mathopen{[\![}}
\def\rcro{\mathclose{]\!]}}

\theoremstyle{definition}
\newtheorem{Def}{Definition}
\newtheorem{Not}[Def]{Notation}

\theoremstyle{plain}
\newtheorem{theo}{Theorem}
\newtheorem{prop}[theo]{Proposition}

\newtheorem{lemme}[theo]{Lemma}
\newtheorem{claim}{Claim}

\theoremstyle{plain}

\theoremstyle{remark}
\newtheorem{Rems}{Remarks}
\newtheorem{Rem}[Rems]{Remark}

\title{Large affine spaces of matrices with rank bounded below}
\author{Cl\'ement de Seguins Pazzis\footnote{Professor of Mathematics at Lyc\'ee Priv\'e Sainte-Genevi\`eve, 2, rue
de l'\'Ecole des Postes, 78029 Versailles Cedex, FRANCE.}
\footnote{e-mail address: dsp.prof@gmail.com}}

\begin{document}

\thispagestyle{plain}

\maketitle
\begin{abstract}
Let $\K$ be an arbitrary (commutative) field with at least three elements,
and let $n$, $p$ and $r$ be positive integers with $r \leq \min(n,p)$. In a recent work \cite{dSPlargerank}, we have proved that
an affine subspace of $\Mat_{n,p}(\K)$ containing only matrices of rank
greater than or equal to $r$ must have a codimension greater than or equal to $\binom{r+1}{2}$.
Here, we classify, up to equivalence, these subspaces of minimal codimension $\binom{r+1}{2}$.
This uses our recent classification \cite{dSPlargeaffine} of the affine subspaces of $\Mat_r(\K)$
contained in $\GL_r(\K)$ which have the maximal dimension $\binom{r}{2}$.
\end{abstract}



\section{Introduction}

In this article, we let $\K$ be an arbitrary (commutative) field. We denote by
$\Mat_n(\K)$ the algebra of square matrices with $n$ rows
and entries in $\K$, and by $\GL_n(\K)$ its group of invertible elements.
We denote by $\Mat_{n,p}(\K)$ the vector space of matrices with $n$ rows, $p$ columns and entries in $\K$.
The rank of a matrix $M$ is denoted by $\rk(M)$, and its transpose is denoted by $M^T$.

Let $E$ be a vector space. An \textbf{affine subspace} $\calV$ of $E$ is the image of a linear subspace
$V$ of $E$ under a translation. In that case, one has $\calV=M+V$ for any $M \in \calV$,
and $V$ is uniquely determined by $\calV$ and is called its \textbf{translation vector space}
(it may be seen as the set of vectors $x \in E$ for which $\calV+x=\calV$).

Given two linear (or affine) subspaces $V$ and $W$ of $\Mat_{n,p}(\K)$, we say that $V$ and $W$ are \textbf{equivalent},
and we write $V \sim W$, if $W=PVQ$ for some $(P,Q)\in \GL_n(\K) \times \GL_p(\K)$.
Two matrices $A$ and $B$ of $\Mat_n(\K)$ are called \textbf{congruent}, and we write
$A \approx B$, if $A=PBP^T$ for some $P \in \GL_n(\K)$. Two quadratic forms $q$ and $q'$ on vector spaces over $\K$ are called \textbf{similar}
when $q'$ is equivalent to $\lambda\,q$ for some $\lambda \in \K \setminus \{0\}$.

\vskip 3mm
Linear subspaces of rectangular matrices with conditions on the rank of their elements have been
extensively studied in the last sixty years. Spaces of matrices with rank bounded above
where first investigated by Dieudonn\'e \cite{Dieudonne} and Flanders \cite{Flanders}.
Flanders showed that if $\# \K \geq r$,  a linear subspace of $\Mat_{n,p}(\K)$, with $n \geq p$, consisting of matrices of rank lesser than or equal to $r$
must have dimension bounded above by $nr$, and equality occurs only for the linear subspaces that are equivalent
to the one of matrices with all last $p-r$ columns equal to zero (or to its transpose in the case $n=p$):
see \cite{Meshulam} for the case of an arbitrary field, and
\cite{dSPaffpres} for an extension to affine subspaces.

Linear subspaces of matrices where all the non-zero elements have rank bounded below by some $r \geq 2$
(or where all the non-zero elements have rank equal to some $r>0$)
have also been under extensive scrutiny: in this case, the results depend greatly on the underlying field.
Tools from algebraic topology are commonly used in those works for the fields of complex and real numbers, see for example
\cite{Causin,Meshulam3,Petrovic,Westwick}. For more general settings, methods from algebraic geometry may be
involved (see the seminal paper \cite{Eisenbud}).

\vskip 2mm
For a non-empty subset $\calX$ of $\Mat_{n,p}(\K)$, we define the \textbf{lower rank} of $\calX$ as:
$$\lrk(\calX):=\min\bigl\{\rk M \mid M \in \calX\bigr\}.$$
Here, we will study the affine subspaces of $\Mat_{n,p}(\K)$ such that $\lrk(\calV) \geq r$,
for a fixed $r \in \lcro 1,\min(r,p)\rcro$, a problem linked to the question of whether a linear subspace of
$\Mat_{n,p}(\K)$ is spanned by its matrices of rank lesser than $r$.
The question first arose in a paper of Meshulam \cite{Meshulam2},
where the following result was proved in the special cases where $\K$ is the field of real numbers or an algebraically closed field.

\begin{theo}\label{majotheorem}
Let $n,p,r$ be positive integers with $r \leq \min(n,p)$.
Let $\calV$ be an affine subspace of $\Mat_{n,p}(\K)$ such that $\lrk(\calV) \geq r$.
Then $\codim \calV \geq \binom{r+1}{2}$.
\end{theo}

For the case of an arbitrary field, see the independent proofs of
R. Quinlan \cite{Quinlan} and the author \cite{dSPlargerank}.

Notice that the lower bound in Theorem \ref{majotheorem} is tight, since equality is obtained with the space of all $n \times p$ matrices of the form
$$\begin{bmatrix}
I_r+U & [?]_{r \times (p-r)} \\
[?]_{(n-r)\times r} & [?]_{(n-r) \times (p-r)}
\end{bmatrix}$$
where $I_r$ is the identity matrix of $\Mat_r(\K)$ and $U$ is an arbitrary strictly upper-triangular matrix of $\Mat_r(\K)$
(and we impose no condition on the blocks represented by question marks).
Remark also that if equality occurs in the above theorem, then $\lrk(\calV)=r$.

Our aim here is to classify, up to equivalence, the affine subspaces of lower rank $r$ and of minimal codimension in $\Mat_{n,p}(\K)$.
Our starting point is our recent classification of the $\binom{n}{2}$-dimensional affine subspaces
of non-singular matrices of $\Mat_n(\K)$ (see \cite{dSPlargeaffine}). We recall a few definitions and notations before stating
our results:

\begin{Def}
An affine subspace of $\Mat_n(\K)$ which is included in $\GL_n(\K)$ is called \textbf{maximal} when its dimension is $\binom{n}{2}$.
\end{Def}

Note that this notion of maximality should not be confused with maximality with respect to the inclusion of
affine subspaces.

\begin{Def}
A non-singular matrix $P \in \GL_n(\K)$ is called \textbf{non-isotropic} when
the quadratic form $X \mapsto X^TPX$ is non-isotropic, i.e., $\forall X \in \K^n \setminus \{0\}, \; X^TPX \neq 0$.
\end{Def}

\begin{Not}
We denote by $\Mata_n(\K)$ the set of alternate matrices in $\Mat_n(\K)$, i.e., of matrices
$A$ such that $\forall X\in \K^n, \; X^TAX=0$, i.e., of skew-symmetric matrices with a zero diagonal.
\end{Not}

Classically, the rank of an alternate matrix is even. In particular, no alternate matrix has rank $1$. 

\begin{Not}
Given respective subsets $V$ and $W$ of $\Mat_n(\K)$ and $\Mat_p(\K)$, we set
$$V \vee W :=\biggl\{\begin{bmatrix}
A & B \\
0 & C
\end{bmatrix} \mid (A,B,C)\in V \times \Mat_{n,p}(\K) \times W
\biggr\} \subset \Mat_{n+p}(\K).$$
\end{Not}

\begin{theo}[Classification theorem for maximal affine subspaces of non-singular matrices]\label{largeaffinenonsingular}
Let $n$ be a positive integer. Assume that $\# \K \geq 3$.
\begin{enumerate}[(a)]
\item Let $(P_1,\dots,P_p)\in \GL_{n_1}(\K) \times \cdots \times \GL_{n_p}(\K)$ be a list of
non-isotropic matrices with $n_1+\cdots+n_p=n$. Then
$$I_n+\Bigl(P_1\Mata_{n_1} \vee \cdots \vee P_p\Mata_{n_p}\Bigr)$$
is a maximal affine subspace of non-singular matrices of $\Mat_n(\K)$.
\item Conversely, let $\calV$ be a maximal affine subspace of non-singular matrices of $\Mat_n(\K)$.
Then there is a list $(P_1,\dots,P_p)\in \GL_{n_1}(\K) \times \cdots \times \GL_{n_p}(\K)$ of
non-isotropic matrices such that
$$\calV \; \sim \; I_n+\Bigl(P_1\Mata_{n_1} \vee \cdots \vee P_p\Mata_{n_p}\Bigr).$$
Moreover, given another list $(Q_1,\dots,Q_q)\in \GL_{m_1}(\K) \times \cdots \times \GL_{m_q}(\K)$ of non-isotropic matrices,
one has
$$\calV \; \sim \; I_n+\Bigl(Q_1\Mata_{m_1} \vee \cdots \vee Q_q\Mata_{m_q}\Bigr)$$
if and only if $p=q$, and, for every $k \in \lcro 1,p\rcro$, the quadratic form $X \mapsto X^TP_kX$ is similar to $X \mapsto X^TQ_k X$.
\end{enumerate}
\end{theo}

The theorem fails for $\# \K=2$ (see Section 6 of \cite{dSPlargeaffine}),
and no classification in this case is known to this day.

\begin{Def}
A maximal affine subspace $\calV$ of $\Mat_n(\K)$ is called \textbf{reducible} when there exists
an integer $p \in \lcro 1,n-1\rcro$, a maximal affine subspace $\calV_1$ of non-singular matrices of $\Mat_p(\K)$
and a maximal affine subspace $\calV_2$ of non-singular matrices of $\Mat_{n-p}(\K)$ such that
$\calV \sim \calV_1 \vee \calV_2$. \\
Otherwise, $\calV$ is called irreducible.
\end{Def}

\noindent As a consequence of Theorem \ref{largeaffinenonsingular}, $\calV$ is irreducible if and only if
$\calV  \sim I_n+P\Mata_n(\K)$ for some non-isotropic matrix $P \in \GL_n(\K)$.

\begin{Not}
Let $n,p,r$ be positive integers, with $n \geq r$ and $p \geq r$.
Given a subset $\calX$ of $\Mat_r(\K)$, we denote by
$i_{n,p}(\calX)$ the set of all matrices of $\Mat_{n,p}(\K)$ of the form
$$\begin{bmatrix}
A & [?]_{r \times (p-r)} \\
[?]_{(n-r) \times r} & [?]_{(n-r) \times (p-r)}
\end{bmatrix} \quad \text{with $A \in \calX$.}$$
\end{Not}

\noindent Here is the theorem we wish to prove:

\begin{theo}[Classification theorem for maximal affine subspaces with rank bounded below]
\label{classaffineboundedbelow}
Let $n$ and $p$ be two positive integers, and $r \in \lcro 2,\min(n,p)\rcro$. Assume that $\# \K \geq 3$.
\begin{enumerate}[(a)]
\item Let $\calW_1$ be a maximal affine subspace of non-singular matrices of $\Mat_r(\K)$.
Then $\calV_1:=i_{n,p}(\calW_1)$ is an affine subspace of $\Mat_{n,p}(\K)$ such that
$\lrk\bigl(\calV_1\bigr)=r$ and $\codim \calV_1=\binom{r+1}{2}$.
\item Conversely, let $\calV$ be an affine subspace of $\Mat_{n,p}(\K)$ such that
$\lrk(\calV)=r$ and $\codim \calV=\binom{r+1}{2}$. Then there exists a maximal affine subspace
$\calW$ of non-singular matrices of $\Mat_r(\K)$ such that
$$\calV \sim i_{n,p}(\calW).$$
Moreover, given another maximal affine subspace
$\calW'$ of non-singular matrices of $\Mat_r(\K)$, one has
$$\calV \sim i_{n,p}(\calW') \; \Longleftrightarrow \; \calW \sim \calW'.$$
\end{enumerate}
\end{theo}

Note that point (a) is an easy observation.

\begin{Rem}
The above theorem fails for $r=1$ and $\min(n,p)>1$, but in this case the classification is easy,
since we are then dealing with the \emph{non-linear} hyperplanes of $\Mat_{n,p}(\K)$.
For every such hyperplane $\calH$, there is a unique matrix $A \in \Mat_{n,p}(\K)$ such that
$\calH=\bigl\{M \in \Mat_{n,p}(\K) : \; \tr(A^TM)=1\bigr\}$, and the equivalence class of $\calH$ is uniquely determined by
that of $A$, i.e., by the rank of $A$.
\end{Rem}

\vskip 2mm
Using Theorem \ref{largeaffinenonsingular}, it follows that the affine subspaces of $\Mat_{n,p}(\K)$ with codimension $\binom{r+1}{2}$ and
lower rank $r \in \lcro 2,\min(n,p)\rcro$ are classified, up to equivalence, by the lists of the form
$\bigl([\varphi_1],\dots,[\varphi_q]\bigr)$, where
$\varphi_1,\dots,\varphi_q$ are non-isotropic quadratic forms over $\K$ with
$\underset{k=1}{\overset{q}{\sum}} \dim \varphi_k=r$, and the $[\varphi_k]$'s denote their similarity classes.

\vskip 2mm
\noindent \textbf{Structure of the proof:}
We will start by proving the uniqueness statement in point (b) of Theorem \ref{classaffineboundedbelow} as it is fairly easy
(Section \ref{uniquenesssection}).
The proof of the existence statement is much more complicated and will involve Theorem \ref{largeaffinenonsingular}:
the basic strategy is laid out in Paragraph \ref{strategy};
the main difficulties lie in the case $r=p=n-1$, which will require an induction and will constitute the main part of the proof (Paragraphs
\ref{r=p=n-1start}, \ref{r=p=n-1irr} and \ref{r=p=n-1generalsection});
from there, the general case will follow rather easily (Paragraph \ref{generalcasesection}).

\section{The uniqueness statement}\label{uniquenesssection}

Let $n$, $p$ and $r$ be three positive integers with $\min(n,p) \geq r$.
Let $\calW$ and $\calW'$ be two maximal affine subspaces of non-singular matrices of $\Mat_r(\K)$. \\
Assume first that $\calW \sim \calW'$. Then $\calW'=P\,\calW \,Q$ for some $(P,Q)\in \GL_r(\K)^2$.
Setting
$$P_1:=P \oplus I_{n-r} \in \GL_n(\K) \quad \text{and} \quad 
Q_1:=Q \oplus I_{p-r} \in \GL_p(\K),$$
we obviously have $i_{n,p}(\calW')=P_1\,i_{n,p}(\calW)\,Q_1$.

\vskip 3mm
Conversely, assume that $i_{n,p}(\calW')=P_1\,i_{n,p}(\calW)\,Q_1$ for some $(P_1,Q_1)\in \GL_n(\K) \times \GL_p(\K)$.
Denote by $(e_1,\dots,e_n)$ the canonical basis of $\K^n$. Denote by $W$ (respectively by $W'$)
the translation vector space of $\calW$ (respectively of $\calW'$). \\
Define $K_W$ as the set of all vectors $x \in \Vect(e_1,\dots,e_n)$ such that
$i_{n,p}(W)$ contains every matrix of $\Mat_{n,p}(\K)$ with column space $\Vect(x)$.

\begin{claim}
One has $K_W=\Vect(e_{r+1},\dots,e_n)$.
\end{claim}

\begin{proof}Choose $A \in \calW$.
Notice that $K_W$ is a linear subspace of $\Vect(e_1,\dots,e_n)$ and that it obviously contains
$\Vect(e_{r+1},\dots,e_n)$. It thus suffices to show that $K_W \cap \Vect(e_1,\dots,e_r)=\{0\}$.
Assume on the contrary that there exists a non-zero vector $y \in K_W \cap \Vect(e_1,\dots,e_r)$.
Considering $y$ as a vector of $\K^r$, we may find a non-zero vector $x$ of $\K^r$ such that
$Ax=y$ (as $A$ is non-singular). Some rank $1$ matrix $B$ of $\Mat_r(\K)$ then satisfies $Bx=y$,
and it follows from the assumption $y \in K_W$ that $i_{n,p}(W)$ contains $\begin{bmatrix}
B & [0]_{r\times (p-r)} \\
[0]_{(n-r) \times r} & [0]_{(n-r) \times (p-r)}
\end{bmatrix}$, showing that $B \in W$. However $A-B \in \calW$ and $A-B$ is singular as $(A-B)x=0$. This is a contradiction.
\end{proof}

Similarly, one has $K_{W'}=\Vect(e_{r+1},\dots,e_n)$.
Now equality $i_{n,p}(\calW')=P_1\,i_{n,p}(\calW)\,Q_1$ yields
$i_{n,p}(W')=P_1\,i_{n,p}(W)\,Q_1$ and $K_{W'}=P_1\,K_{W}$, which shows that $P_1=\begin{bmatrix}
P & [0]_{r \times (n-r)} \\
[?]_{(n-r) \times r} & P'
\end{bmatrix}$ for some $(P,P')\in \GL_r(\K) \times \GL_{n-r}(\K)$. \\
On the other hand, transposing the above equality yields $i_{p,n}((\calW')^T)=Q_1^T\,i_{p,n}(\calW^T)\,P_1^T$
and applying the previous result yields that
$Q_1^T=\begin{bmatrix}
Q & [0]_{r \times (p-r)} \\
[?]_{(p-r) \times r} & Q'
\end{bmatrix}$ for some $(Q,Q')\in \GL_r(\K) \times \GL_{p-r}(\K)$.
Then
$$P_1\,i_{n,p}(\calW)\,Q_1=i_{n,p}(P\,\calW\,Q^T)$$
and it immediately follows that $\calW'=P\,\calW\,Q^T$. Therefore, $\calW \sim \calW'$, QED.

\section{The existence statement}\label{existencesection}

In the whole section, we let $n$, $p$ and $r$ be three positive integers with $\min(n,p)\geq r$,
and we assume $\# \K \geq 3$.
Let $\calV$ be an affine subspace of $\Mat_{n,p}(\K)$ such that $\lrk(\calV)=r$ and $\codim \calV=\binom{r+1}{2}$. \\
We wish to find a maximal affine subspace $\calW$ of non-singular matrices of $\Mat_r(\K)$ such that
$\calV \sim i_{n,p}(\calW)$, i.e., $\calV=P\,i_{n,p}(\calW)\,Q$
for some pair $(P,Q) \in \GL_n(\K) \times \GL_p(\K)$. \\
This will be achieved using a series of steps involving reductions of the following (essentially equivalent) types, which all transform
$\calV$ into an equivalent affine subspace:
\begin{itemize}
\item left and right-multiplication of $\calV$ with non-singular square matrices;
\item row and column operations.
\end{itemize}

\subsection{Putting $\calV$ to a roughly-reduced form}\label{strategy}

Here, we follow some ideas of Meshulam \cite{Meshulam2}.
Denote by $V$ the translation vector space of $\calV$. 
We say that $\calV$ is \textbf{roughly-reduced} when it contains the matrix
$$J:=\begin{bmatrix}
I_r & [0]_{r \times (p-r)} \\
[0]_{(n-r) \times r} & [0]_{(n-r) \times (p-r)}
\end{bmatrix}.$$

Since $\calV$ contains a rank $r$ matrix, which is equivalent to $J$,
we may find a pair $(P,Q)\in \GL_n(\K) \times \GL_p(\K)$ such that $P\calV Q$ contains $J$.
Therefore, in the rest of the proof, we lose no generality in assuming that $\calV$ is roughly-reduced.

\noindent In that case, we denote by $W$ the linear subspace of the matrices $M \in \Mat_r(\K)$ such that
$\begin{bmatrix}
M & [0]_{r \times (p-r)} \\
[0]_{(n-r) \times r} & [0]_{(n-r) \times (p-r)}
\end{bmatrix}$ belongs to $V$.
Note that $\calV$ contains $\begin{bmatrix}
I_r+M & [0]_{r \times (p-r)} \\
[0]_{(n-r) \times r} & [0]_{(n-r) \times (p-r)}
\end{bmatrix}$ for every $M \in W$, hence
$$\calW:=I_r+W$$
is an affine subspace of non-singular matrices of $\Mat_r(\K)$: we will call $\calW$ the \textbf{core space} of $\calV$.

\vskip 2mm
We now define $H_\calV$ as the linear subspace of all triples $(B,C,D)\in \Mat_{n-r,r}(\K) \times \Mat_{r,p-r}(\K) \times \Mat_{n-r,p-r}(\K)$
such that $V$ contains a matrix of the form
$\begin{bmatrix}
[?]_{r \times r} & C \\
B & D
\end{bmatrix}$.
The rank theorem then shows that
$$\dim \calV=\dim W+\dim H_\calV=\dim \calW+\dim H_\calV.$$
However $\dim \calW \leq \binom{r}{2}$ by Theorem \ref{majotheorem}, whilst
$\dim H_\calV \leq np-r^2$ and
$\binom{r}{2}+np-r^2=np-\binom{r+1}{2}=\dim \calV$,
hence
$$\dim \calW=\binom{r}{2} \quad \text{and} \quad \dim H_\calV=np-r^2.$$
It follows that:
\begin{enumerate}[(i)]
\item $\calW$ is a maximal affine subspace of non-singular matrices of $\Mat_r(\K)$.
\item For every $(B,C,D)\in \Mat_{n-r,r}(\K) \times \Mat_{r,p-r}(\K) \times \Mat_{n-r,p-r}(\K)$,
the linear subspace $V$ contains a matrix of the form $\begin{bmatrix}
[?]_{r \times r} & C \\
B & D
\end{bmatrix}$.
\end{enumerate}

From there, we aim at proving that $\calV \sim i_{n,p}(\calW)$.
We will first do this in the case $p=r$ and $n=p+1$ (see the next three sections),
and then generalize the result in the case $r \geq 2$.

\subsection{The existence statement for $r=p=n-1$ (I): general considerations}\label{r=p=n-1start}

Proving the following proposition will obviously solve the existence problem in the case
$r=p=n-1$, using the considerations from Paragraph \ref{strategy}.

\begin{prop}\label{caser=p=n-1}
Let $\calV$ be an affine subspace of $\Mat_{n+1,n}(\K)$ in which every matrix has rank $n$
and such that $\codim \calV=\binom{n+1}{2}$. \\
Assume that $\calV$ is roughly-reduced and denote by $\calW$ its core space. \\
Then there is a list $(\lambda_1,\dots,\lambda_n)\in \K^n$ such that the series of row operations
$L_i \leftarrow L_i+\lambda_i\,L_{n+1}$, for $i$ from $1$ to $n$,
transforms $\calV$ into $i_{n+1,n}(\calW)$.
\end{prop}

In the proof, we denote by $W$ the translation vector space of $\calW$ (as in the previous section).
Our proof of Proposition \ref{caser=p=n-1} will be done by induction on $n$, with two stages:
we will start with the case where $W$ is irreducible as a maximal affine subspace of non-singular matrices of $\Mat_n(\K)$,
and we will then move on to the case where $W$ is reducible.

Before giving the details of the proof, let us recall the following basic results
(see \cite{dSPlargeaffine} for the proofs):
\begin{enumerate}[(i)]
\item For every $(P,Q)\in \GL_n(\K)^2$ and every $\lambda \in \K \setminus \{0\}$,
one has $P\Mata_n(\K)Q^{-1}=\lambda (PQ^T) \Mata_n(\K)$.
\item For every $X \in \K^n \setminus \{0\}$, one has $\Mata_n(\K) X=\{X\}^\bot$, where orthogonality has to be understood
for the non-degenerate symmetric bilinear form $(Y,Z) \mapsto Y^TZ$.
\end{enumerate}

\begin{Rem}\label{reductionremark}
In the course of the proof, it will be necessary to modify $\calV$ so that
its core space has a reduced shape. Let us see how to justify this. Let $(Q_1,Q_2) \in \GL_n(\K)^2$ and set
$\widetilde{Q_1}:=Q_1 \oplus 1 \in \GL_{n+1}(\K)$. Assume that $Q_1 \calW Q_2$ contains $I_n$. Then, with the assumptions of Proposition \ref{caser=p=n-1},
the affine space $\calV^{(1)}:=\widetilde{Q_1}\calV Q_2$ is still roughly-reduced and its core space is
$\calW^{(1)}:=Q_1 \calW Q_2$. Assume that there exists a list $(\lambda_1,\dots,\lambda_n)\in \K^n$ such that the series of row operations
$L_i \leftarrow L_i+\lambda_i\,L_{n+1}$, for $i$ from $1$ to $n$,
transforms $\calV^{(1)}$ into $i_{n+1,n}(\calW^{(1)})$. Setting $C:=\begin{bmatrix}
\lambda_1 \\
\vdots \\
\lambda_n
\end{bmatrix}$, this means that
$R\,\calV^{(1)}=i_{n+1,n}(\calW^{(1)})$ for $R:=\begin{bmatrix}
I_n & C \\
[0]_{1 \times n} & 1
\end{bmatrix}$. Then $R':=\widetilde{Q_1}^{-1}\,R\,\widetilde{Q_1}=\begin{bmatrix}
I_n & Q_1^{-1} C \\
[0]_{1 \times n} & 1
\end{bmatrix}$ satisfies $R'\,\calV=i_{n+1,n}(\calW)$, which means, if we write $Q_1^{-1}C=\begin{bmatrix}
\lambda'_1 \\
\vdots \\
\lambda'_n
\end{bmatrix}$, that the series of row operations
$L_i \leftarrow L_i+\lambda'_i\,L_{n+1}$, for $i$ from $1$ to $n$,
transforms $\calV$ into $i_{n+1,n}(\calW)$.

We conclude that, in proving Proposition \ref{caser=p=n-1}, no generality is lost by replacing
the core space $\calW$ of $\calV$ with any equivalent affine subspace of $\Mat_n(\K)$ containing $I_n$.
With the above explanation, this may be obtained by taking $Q_1=Q_2^{-1}$; in this case,
if in addition $\calW=I_n+P\Mata_n(\K)$ for some non-isotropic matrix $P$, then, for every $\lambda \in \K \setminus \{0\}$, one has
$Q_1 \calW Q_1^{-1}=I_n+(\lambda Q_1PQ_1^T)\Mata_n(\K)$, and $\lambda Q_1PQ_1^T$ is still non-isotropic.
\end{Rem}

\subsection{The existence statement for $r=p=n-1$ (II): the irreducible case}\label{r=p=n-1irr}

In the case where $\calW$ is irreducible, Proposition \ref{caser=p=n-1} may be restated as follows
(see Remark \ref{reductionremark} of Paragraph \ref{r=p=n-1start}).

\begin{prop}\label{irrprop}
Let $P \in \GL_n(\K)$ be a non-isotropic matrix and set $\calW:=I_n+P\,\Mata_n(\K)$.
Let $\calV$ be an affine subspace of $\Mat_{n+1,n}(\K)$ with codimension $\binom{n+1}{2}$ in which every matrix has rank $n$,
and which contains the affine subspace
$\biggl\{\begin{bmatrix}
M \\
0
\end{bmatrix} \mid M \in \calW\biggr\}$. Then there is a list $(\lambda_1,\dots,\lambda_n)\in \K^n$ such that the series of row operations
$L_i \leftarrow L_i+\lambda_i\,L_{n+1}$, for $i$ from $1$ to $n$,
transforms $\calV$ into $i_{n+1,n}(\calW)$.
\end{prop}

Before proving this proposition, we need a lemma:

\begin{lemme}\label{rank1lemma}
Let $\calV$ satisfying the assumptions of Proposition \ref{caser=p=n-1}
and assume that the translation vector space $W$ of its core space $\calW$
has the form $P\Mata_n(\K)$ for some $P \in \GL_n(\K)$. \\
Then, for every $L \in \Mat_{1,n}(\K) \setminus \{0\}$,
there is a rank $1$ matrix with $L$ as last row in the translation vector space $V$ of $\calV$ .
\end{lemme}

\begin{proof}
Let $L \in \Mat_{1,n}(\K) \setminus \{0\}$.
There is a non-singular matrix $Q \in \GL_n(\K)$ such that
$$L\,Q=L_1:=\begin{bmatrix}
1 & 0 & \cdots & 0
\end{bmatrix} \in \Mat_{1,n}(\K).$$
Note that $\calV Q$ still has codimension $\binom{n+1}{2}$ in $\Mat_{n+1,n}(\K)$, has lower rank $n$, is
roughly-reduced, and the translation vector space of its core space
is $P\Mata_n(\K)Q=\bigl(P(Q^T)^{-1}\bigr)\,\Mata_n(\K)$.
If $V Q$ contains a rank $1$ matrix $A$ with $L_1$ as last row, then $V$ contains the rank $1$ matrix $AQ^{-1}$ with
$L_1 Q^{-1}=L$ as last row. \\
Therefore, it suffices to tackle the case $L=L_1$, which we now do. \\
By point (ii) of Paragraph \ref{strategy}, we know that $V$ contains a matrix of the form
$$\begin{bmatrix}
[?]_{n \times 1} & N_1 \\
1 & [0]_{1 \times (n-1)}
\end{bmatrix} \quad \text{with $N_1 \in \Mat_{n,n-1}(\K)$.}$$
Write every $M \in \Mat_n(\K)$ as
$M=\begin{bmatrix}
[?]_{n \times 1} & R(M) \\
\end{bmatrix}$ with $R(M) \in \Mat_{n,n-1}(\K)$ and note that $M \mapsto R(M)$ is one-to-one on $\calW$ since no matrix of $\Mata_n(\K)$ has rank $1$
and $P$ is non-singular.
For every $M \in \calW$ and every $\lambda \in \K$,
the affine subspace $\calV$ contains a matrix of the form
$\begin{bmatrix}
[?]_{n \times 1} & R(M)+\lambda N_1 \\
\lambda & [0]_{1 \times (n-1)}
\end{bmatrix}$, which shows that $R(M)+\lambda N_1$ has rank $n-1$.
However, $R(\calW)$ is an affine subspace of $\Mat_{n,n-1}(\K)$ with lower rank $n-1$ and
$\dim R(\calW)=\dim \calW=(n-1)+\binom{n-1}{2}$. Theorem \ref{majotheorem}
thus shows that $R(\calW)$ is maximal among the affine subspaces of $\Mat_{n,n-1}(\K)$ with lower rank $n-1$. It follows that $N_1$ belongs to $R(W)$,
hence $V$ contains a matrix of the form
$\begin{bmatrix}
[?]_{n \times 1} & [0]_{n \times (n-1)} \\
1 & [0]_{1 \times (n-1)}
\end{bmatrix} \in \Mat_{n+1,n}(\K)$, QED.
\end{proof}

We are now ready to prove Proposition \ref{irrprop} by induction on $n$.
If $n=1$, then the translation vector space $V$ of $\calV$ is a $1$-dimensional subspace spanned by a matrix
of the form $\begin{bmatrix}
a \\
1
\end{bmatrix}$, whereas $\calW=\{I_1\}$: it thus suffices to use the row operation
$L_1 \leftarrow L_1-aL_2$. \\
Let $n \geq 2$ and assume that the result of Proposition \ref{irrprop} holds in $\Mat_{n,n-1}(\K)$.
Let $\calV \subset \Mat_{n+1,n}(\K)$, $\calW \subset \Mat_n(\K)$ and $P \in \GL_n(\K)$ satisfying the assumptions of Proposition \ref{irrprop}.
Denote by $V$ the translation vector space of $\calV$.

By Remark \ref{reductionremark} of Paragraph \ref{r=p=n-1start},
we lose no generality in replacing $P$ with $\lambda\,QPQ^T$ for an arbitrary $Q \in \GL_n(\K)$ and a non-zero scalar $\lambda$.
Since $X \mapsto X^TPX$ is non-isotropic, we may complete the first vector $e_1$ of the canonical basis of $\K^n$
to a basis $(e_1,f_2,\dots,f_n)$ such that $(e_1)^TPf_i=0$ for every $i \in \lcro 2,n\rcro$.
This shows that we lose no generality in assuming that
$$P=\begin{bmatrix}
1 & [0]_{1 \times (n-1)} \\
C & Z
\end{bmatrix} \quad \text{for some $(C,Z)\in \Mat_{n-1,1}(\K) \times \GL_{n-1}(\K)$,}$$
where $Z$ must be non-isotropic.
In this case, we have:

\begin{center}
(R1) : The matrices of $V$ with zero as $(n+1)$-th row are precisely the matrices
$$\begin{bmatrix}
0 & L \\
-ZL^T & CL+ZA \\
0 & [0]_{1 \times (n-1)}
\end{bmatrix} \quad \text{with $(L,A)\in \Mat_{1,n-1}(\K) \times \Mata_{n-1}(\K)$.}$$
\end{center}

\vskip 5mm
From there, our short-term aim is to find appropriate scalars $\lambda_1,\dots,\lambda_n$
for the row operations mentioned in the statement of Proposition \ref{irrprop}.
Note that performing any row operation of the form $L_i \leftarrow L_i+\lambda L_{n+1}$ (where $i \neq n+1$)
leaves the core space of $\calV$ unchanged.
We start by looking for appropriate values for $\lambda_2,\dots,\lambda_n$.

\vskip 5mm
For $M \in \calV$, write now
$M=\begin{bmatrix}
[?]_{(n+1) \times 1} & K(M)
\end{bmatrix}$ with $K(M) \in \Mat_{n+1,n-1}(\K)$. Since
every matrix of $\calV$ has rank $n$, its columns are linearly independent,
hence $\calV':=K(\calV)$ is an affine subspace of $\Mat_{n+1,n-1}(\K)$ in which every matrix has rank $n-1$.
Notice that, for every $(L,A)\in \Mat_{1,n-1}(\K) \times \Mata_{n-1}(\K)$, the affine subspace $\calV'$ contains
$$\begin{bmatrix}
L \\
I_{n-1}+CL+ZA \\
[0]_{1 \times (n-1)}
\end{bmatrix} \in \Mat_{n+1,n-1}(\K)$$
and its translation vector space $V'$ contains the matrix
$$\begin{bmatrix}
L \\
CL+ZA \\
[0]_{1 \times (n-1)}
\end{bmatrix}.$$
Applying Lemma \ref{rank1lemma} to $\calV$ and adding a well-chosen matrix of $V'$,
we deduce that, for every $L \in \Mat_{1,n-1}(\K)$, the vector space $V'$
contains a matrix of the form
$$\begin{bmatrix}
[0]_{1\times (n-1)} \\
[?]_{(n-1)\times (n-1)} \\
L
\end{bmatrix}.$$ \\
Now, consider the affine subspace $\calV_1$ of $\Mat_{n,n-1}(\K)$ consisting of the matrices of $\calV'$ with zero as first row:
write every such matrix as $M=\begin{bmatrix}
[0]_{1 \times (n-1)} \\
J(M)
\end{bmatrix}$ with $J(M) \in \Mat_{n,n-1}(\K)$. Note that the affine subspace $J(\calV_1)$ satisfies all the conditions
of Proposition \ref{irrprop}, with $I_{n-1}+Z\Mata_{n-1}(\K)$ as core space. \\
Applying the induction hypothesis, we see that, by using a series of row operations on $\calV$ of the form $L_i \leftarrow L_i+\lambda_i L_{n+1}$,
for $i \in \lcro 2,n+1\rcro$, the situation is reduced to the one where
\begin{equation}\label{structuredeV1}
\calV_1=\Biggl\{\begin{bmatrix}
[0]_{1 \times (n-1)} \\
I_{n-1}+ZA \\
L
\end{bmatrix} \mid (L,A)\in \Mat_{1,n-1}(\K) \times \Mata_{n-1}(\K)\Biggr\},
\end{equation}
with $\calW$ left unchanged.

\vskip 5mm
We now search for an appropriate $\lambda_1$.
Let $L \in \Mat_{1,n-1}(\K)$. We know from Lemma \ref{rank1lemma} (the case $L=0$ being trivial) that
we may find some $\alpha(L) \in \K$ such that $V$ contains a matrix of the form
$$\begin{bmatrix}
0 & \alpha(L)\cdot L \\
[0]_{(n-1) \times 1} & [?]_{(n-1)\times (n-1)} \\
0 & L
\end{bmatrix}.$$

\begin{claim}
The map $L \mapsto \alpha(L)$ is constant on $\Mat_{1,n-1}(\K) \setminus \{0\}$.
\end{claim}

\begin{proof}
Let $L \in \Mat_{1,n-1}(\K)$. By summing a matrix in $V$ of the form $\begin{bmatrix}
0 & \alpha(L)\!\cdot\!L \\
[0]_{(n-1) \times 1} & [?]_{(n-1)\times (n-1)} \\
0 & L
\end{bmatrix}$ with the matrix $-\alpha(L)\!\cdot\!\begin{bmatrix}
0 & L \\
-ZL^T & CL \\
0 & [0]_{1 \times (n-1)}
\end{bmatrix}$, we obtain a matrix in $V$ which has the form
$\begin{bmatrix}
0 & [0]_{1 \times (n-1)} \\
\alpha(L)\!\cdot\!ZL^T & [?]_{(n-1) \times (n-1)} \\
0 & L
\end{bmatrix}$. Using \eqref{structuredeV1}, we see that this matrix may be written as
$\begin{bmatrix}
0 & [0]_{1 \times (n-1)} \\
\alpha(L)\!\cdot\!ZL^T & ZA \\
0 & L
\end{bmatrix}$ for some $A \in \Mata_{n-1}(\K)$. Summing it with $\begin{bmatrix}
0 & [0]_{1 \times (n-1)} \\
[0]_{(n-1)\times 1} & -ZA \\
0 & [0]_{1 \times (n-1)}
\end{bmatrix} \in V$,
we deduce that
$$\begin{bmatrix}
0 & [0]_{1 \times (n-1)} \\
\alpha(L)\!\cdot\! ZL^T & [0]_{(n-1) \times (n-1)} \\
0 & L
\end{bmatrix}\in V.$$
However, we know from property (R1) that any matrix of $V$ which has zero as first and $(n+1)$-th row
also has zero as first column. It follows that, on the linear subspace of $V$ consisting of its matrices with zero as first row,
the first column is a linear function of the $(n+1)$-th row, which yields that
$L \in \Mat_{1,n-1}(\K) \mapsto \alpha(L)\!\cdot\!ZL^T$ is linear.
Therefore $L \mapsto \alpha(L)\!\cdot\!L$ is an endomorphism of $\Mat_{1,n-1}(\K)$;
as every non-zero vector of $\Mat_{1,n-1}(\K)$ is an eigenvector of it, the map
$L \mapsto \alpha(L)$ is constant on $\Mat_{1,n-1}(\K) \setminus \{0\}$.
\end{proof}

Denote by $\mu$ the sole value of $L \mapsto \alpha(L)$ on $\Mat_{1,n-1}(\K) \setminus \{0\}$.
Performing the row operation $L_1 \leftarrow L_1-\mu\,L_{n+1}$ on $\calV$,
we now have:
\begin{center}
(R2) : $V$ contains a matrix of the form $\begin{bmatrix}
0 & [0]_{1 \times (n-1)} \\
[0]_{(n-1) \times 1} & [?]_{(n-1) \times (n-1)} \\
0 & L
\end{bmatrix}$ for every $L \in \Mat_{1,n-1}(\K)$.
\end{center}
Alas, by doing so, we may have lost property \eqref{structuredeV1} for the new $\calV_1$ space!
Nevertheless, we may use an additional series of row operations of the form $L_i \leftarrow L_i+\lambda L_{n+1}$, with
$i \in \lcro 2,n+1\rcro$, so as to recover property \eqref{structuredeV1}.
By performing those row operations, we keep property (R2), so we may now assume that both
\eqref{structuredeV1} and (R2) are satisfied.

Let then $L \in \Mat_{1,n-1}(\K)$. Using (R2) and \eqref{structuredeV1},
we know that $V$ contains $\begin{bmatrix}
0 & [0]_{1 \times (n-1)} \\
[0]_{(n-1) \times 1} & ZA \\
0 & L
\end{bmatrix}$ for some $A \in \Mata_{n-1}(\K)$. Subtracting it with the matrix
$\begin{bmatrix}
0 & [0]_{1 \times (n-1)} \\
[0]_{(n-1) \times 1} & ZA \\
0 & [0]_{1 \times (n-1)}
\end{bmatrix} \in V$, we deduce:

\begin{center}
(R3) : $V$ contains
$\begin{bmatrix}
[0]_{n \times 1} & [0]_{n \times (n-1)}  \\
0 & L
\end{bmatrix}$ for every $L \in \Mat_{1,n-1}(\K)$.
\end{center}
Set now $L_1:=\begin{bmatrix}
1 & 0 & \cdots & 0
\end{bmatrix} \in \Mat_{1,n}(\K)$.
Then we know from Lemma \ref{rank1lemma} that $V$ contains
$\begin{bmatrix}
M_1 \\
L_1
\end{bmatrix}$ for some $M_1 \in \Mat_n(\K)$.

\begin{claim}\label{M1in}
One has $M_1 \in P\Mata_n(\K)$.
\end{claim}

\begin{proof}
Note that, for every $M \in P\Mata_n(\K)$, every $\beta \in \K \setminus \{0\}$ and every $L \in \Mat_{1,n-1}(\K)$,
the affine subspace $\calV$ contains the matrix
$\begin{bmatrix}
I_n+M+\beta M_1 \\
L_2
\end{bmatrix} \in \Mat_{n+1,n}(\K)$ with $L_2:=\begin{bmatrix}
\beta & L
\end{bmatrix}$.
Denote by $(e_1,\dots,e_n)$ the canonical basis of $\Mat_{1,n}(\K)$.
If $I_n+M+\beta M_1$ were singular for some $M \in P\Mata_n(\K)$ and some $\beta \in \K \setminus \{0\}$,
then it would have rank $n-1$: if in addition its first column were non-zero, then
its row space would be a linear hyperplane $H$ of $\Mat_{1,n}(\K)$ which is different from
$\Vect(e_2,\dots,e_n)$, and it would then have a common point with the non-parallel
affine hyperplane $\beta e_1+\Vect(e_2,\dots,e_n)$, i.e., we would be able to find $L \in \Mat_{1,n-1}(\K)$ such that the row matrix
$\begin{bmatrix}
\beta & L
\end{bmatrix}$ belongs to the row space of $I_n+M+\beta M_1$:
this would provide a matrix in $\calV$ with rank $n-1$, a contradiction.

We deduce that, for every $M \in P\Mata_n(\K)$ and every $\beta \in \K \setminus \{0\}$ such that $I_n+M+\beta M_1$ is singular,
the kernel of $I_n+M+\beta M_1$ is spanned by the column matrix $X_1:=\begin{bmatrix}
1 & 0 & \cdots & 0
\end{bmatrix}^T$.

Assume now that $M_1 \not\in P\Mata_n(\K)$,
and set $A_1:=P^{-1}M_1$. The quadratic form $q : X \mapsto X^TA_1X$ on $\K^n$ is then non-zero,
hence we may choose a vector $X \in \K^n \setminus \Vect(X_1)$ for which $q(X) \neq 0$
(assume this is not possible: then on the one hand $q(X_1) \neq 0$; on the other hand,
choosing $X_2 \in \K^n \setminus \Vect(X_1)$, we would find that the quadratic form $q$ vanishes on
three\footnote{Remember that $\# \K \geq 3$.} distinct $1$-dimensional linear subspaces of $\Vect(X_1,X_2)$, hence $q(X_1)=0$, a contradiction).

Then $A_1X \not\in \{X\}^\bot$, and since $\{X\}^\bot=\Mata_n(\K)X$ has codimension $1$ in $\K^n$,
it would follow that $-P^{-1}X=\beta A_1X+NX$ for some $\beta \in \K$ and some $N \in \Mata_n(\K)$.
However $\beta \neq 0$ since $P^{-1}$ is non-isotropic: it would follow that
$I_n+P\,N+\beta\,M_1$ is singular with $X$ in its kernel, which contradicts the above
proof and the choice of $X$. \\
This \emph{reductio ad absurdum} shows that $M_1 \in P\Mata_n(\K)$.
\end{proof}

Using Claim \ref{M1in} and equality $\calW=I_n+P\Mata_n(\K)$, we find that $V$ contains $\begin{bmatrix}
[0]_{n \times n} \\
L_1
\end{bmatrix}$. Combining this with (R3), we deduce that $V$ contains the matrix
$\begin{bmatrix}
[0]_{n \times n} \\
L
\end{bmatrix}$ for every $L \in \Mat_{1,n}(\K)$. It follows that $i_{n+1,n}(\calW) \subset \calV$, and since the dimensions
are equal, we conclude that $i_{n+1,n}(\calW)=\calV$. This completes the proof of Proposition \ref{irrprop}.

\vskip 2mm
Now, the case where $W$ is irreducible is done.

\subsection{The existence statement for $p=r=n-1$ (III): the general case}\label{r=p=n-1generalsection}

We now prove Proposition \ref{caser=p=n-1} in the general case. Again, we use an induction on $n$.
The proof is straightforward for $n=1$ (see the previous section).
Let $n \geq 2$ be an integer, assume that Proposition \ref{caser=p=n-1} holds for any positive
integer lesser than $n$, and let $\calV  \subset \Mat_{n+1,n}(\K)$ and $\calW$ be as in Proposition \ref{caser=p=n-1}.
Denote respectively by $V$ and $W$ the translation vector spaces of $\calV$ and $\calW$.
By Theorem \ref{largeaffinenonsingular}, there is a list of non-isotropic matrices $(P_1,\dots,P_q)\in \GL_{n_1}(\K) \times \cdots \times \GL_{n_q}(\K)$
such that
$$\calW \sim I_n+\bigl(P_1\Mata_{n_1}(\K) \vee \cdots \vee P_q \Mata_{n_q}(\K)\bigr).$$
By Remark \ref{reductionremark}, we lose no generality in assuming:
\begin{center}
(S1) : \quad $\calW=I_n+\bigl(P_1\Mata_{n_1}(\K) \vee \cdots \vee P_q \Mata_{n_q}(\K)\bigr)$.
\end{center}
For convenience, we now set
$$s:=n_1 \quad ,\quad P:=P_1 \quad \text{and} \quad U:=P_2\Mata_{n_2}(\K) \vee \cdots \vee P_q \Mata_{n_q}(\K).$$
The case $s=n$ having been dealt with in Proposition \ref{irrprop}, we now assume that $s<n$. \\
Recall from point (ii) of Paragraph \ref{r=p=n-1start} that, for any row matrix $L \in \Mat_{1,n}(\K)$,
the subspace $V$ contains a matrix with $L$ as last row.

\begin{Not}
We denote by $\calV_1$ the affine subspace of $\calV$ consisting of its matrices of the form
$\begin{bmatrix}
[?]_{n \times s} & [?]_{n \times (n-s)} \\
[?]_{1 \times s} & [0]_{1 \times (n-s)}
\end{bmatrix}$.
\end{Not}

\noindent For $M \in \calV_1$, we write
$$M=\begin{bmatrix}
J(M) & [?]_{s \times (n-s)} \\
[?]_{(n-s) \times s} & K(M) \\
L(M) & [0]_{1 \times (n-s)}
\end{bmatrix} \quad \text{with $(J(M),K(M),L(M))\in \Mat_s(\K) \times \Mat_{n-s}(\K) \times \Mat_{1,s}(\K)$}$$
and we set
$$T(M)=\begin{bmatrix}
J(M) \\
L(M)
\end{bmatrix} \in \Mat_{s+1,s}(\K).$$
Adding
$\begin{bmatrix}
[0]_{s \times s} & B \\
[0]_{(n-s+1) \times s} & [0]_{(n-s+1) \times (n-s)}
\end{bmatrix}$ (which belongs to $V$ by (S1)) for a well-chosen $B \in \Mat_{s,n-s}(\K)$, we deduce, for every $M \in \calV_1$,
that $\calV$ contains a matrix of the form
$$\begin{bmatrix}
J(M) & [0]_{s \times (n-s)} \\
[?]_{(n-s) \times s} & K(M) \\
L(M) & [0]_{1 \times (n-s)}
\end{bmatrix},$$
which successively shows that $\rk K(M)=n-s$ and $\rk T(M)=s$.  \\
Therefore $T(\calV_1)$ is an affine subspace of rank $s$ matrices of $\Mat_{s+1,s}(\K)$.
Notice also that $T(\calV_1)$ contains the matrix
$\begin{bmatrix}
I_s+N \\
[0]_{1 \times s}
\end{bmatrix}$ for every $N \in P\Mata_s(\K)$ and that it contains
a matrix of the form
$\begin{bmatrix}
[?]_{s \times s} \\
L
\end{bmatrix}$ for every $L \in \Mat_{1,s}(\K)$.
It follows that $\codim T(\calV_1) \geq \binom{s+1}{2}$, therefore $\codim T(\calV_1)=\binom{s+1}{2}$ by Theorem \ref{majotheorem}.
We may then apply the induction hypothesis to $T(\calV_1)$:
we deduce, after a series of row operations on $\calV$ of the form $L_i \leftarrow L_i+\lambda_i\,L_{n+1}$
for $i$ from $1$ to $s$, that we lose no generality in assuming:

\begin{center}
(S2) : \quad $T(\calV_1)$ is the set of all matrices of the form
$\begin{bmatrix}
I_s+P\,A \\
L
\end{bmatrix}$ with $A \in \Mata_s(\K)$ and $L \in \Mat_{1,s}(\K)$.
\end{center}

\begin{Not}
For $M \in \calV$, we denote by $R(M) \in \Mat_{n-s+1,n-s}(\K)$ the matrix such that
$$M=\begin{bmatrix}
[?]_{s \times s} & [?]_{s \times (n-s)} \\
[?]_{(n+1-s) \times s} & R(M)
\end{bmatrix}.$$
\end{Not}
By (S2), we know that $V$ contains the matrix
$\begin{bmatrix}
[0]_{s \times s} & B \\
[0]_{(n+1-s) \times s} & [0]_{(n+1-s) \times (n-s)}
\end{bmatrix}$ for every $B \in \Mat_{s,n-s}(\K)$.
It follows that, for every $M \in \calV$, the affine subspace $\calV$ contains a matrix of the form
$\begin{bmatrix}
[?]_{s \times s} &  [0]_{s \times (n-s)} \\
[?]_{(n+1-s) \times s} & R(M)
\end{bmatrix}$, which entails that $\rk R(M)=n-s$. \\
Therefore:
\begin{enumerate}[(i)]
\item $R(\calV)$ is an affine subspace of $\Mat_{n-s+1,n-s}(\K)$ with lower rank $n-s$.
\item For every matrix $N$ of $I_{n-s}+U$, the subspace $R(\calV)$ contains $\begin{bmatrix}
N \\
[0]_{1 \times (n-s)}
\end{bmatrix}$ (see (S1)).
\item For every $L \in \Mat_{1,n-s}(\K)$, the subspace $R(\calV)$ contains a matrix
of the form $\begin{bmatrix}
[?]_{(n-s) \times (n-s)} \\
L
\end{bmatrix}$ (use point (ii) of Paragraph \ref{strategy}).
\end{enumerate}
As before, we may then apply the induction hypothesis to
$R(\calV)$, and deduce that, after a series of row operations on $\calV$ of the form $L_i \leftarrow L_i+\lambda_i\,L_{n+1}$
for $i$ from $s+1$ to $n$, one may assume:
\begin{center}
(S3) : \quad
$R(\calV)$ is the set of all matrices of the form
$\begin{bmatrix}
I_{n-s}+N \\
L
\end{bmatrix}$ with $N \in U$ and $L \in \Mat_{1,n-s}(\K)$.
\end{center}
Note that properties (S1) and (S2) are preserved by performing these row operations.

\vskip 2mm
\noindent Let us now sum the situation up:
\begin{itemize}
\item[(i)] For every $M \in W=P\Mata_s(\K) \vee U$,
the linear subspace $V$ contains the matrix
$\begin{bmatrix}
M \\
[0]_{1 \times n}
\end{bmatrix}$ (see (S1)).
\end{itemize}
Let $L \in \Mat_{1,s}(\K)$. We know from (S2) that $V$ contains 
a matrix of the form 
$\begin{bmatrix}
[0]_{s \times s} & [0]_{s \times (n-s)} \\
[?]_{(n-s) \times s} & N \\
L & [0]_{1 \times (n-s)}
\end{bmatrix}$
with $N \in \Mat_{n-s}(\K)$. Using (S3), we find that $N \in U$. Using point (i) above, 
we know that $V$ contains the matrix 
$\begin{bmatrix}
[0]_{s \times s} & [0]_{s \times (n-s)} \\
[?]_{(n-s) \times s} & -N \\
[0]_{1 \times s} & [0]_{1 \times (n-s)}
\end{bmatrix}$. Adding those two matrices yields that $V$ contains a matrix of the form
$\begin{bmatrix}
[0]_{s \times s} & [0]_{s \times (n-s)} \\
[?]_{(n-s) \times s} & [0]_{(n-s) \times (n-s)} \\
L & [0]_{1 \times (n-s)}
\end{bmatrix}$. 
Therefore:
\begin{itemize}
\item[(ii)] For every $L \in \Mat_{1,s}(\K)$, the linear subspace $V$ contains a matrix of the form
$\begin{bmatrix}
[0]_{s \times s} & [0]_{s \times (n-s)} \\
[?]_{(n-s) \times s} & [0]_{(n-s) \times (n-s)} \\
L & [0]_{1 \times (n-s)}
\end{bmatrix}$. 
\end{itemize}
Let finally $L' \in \Mat_{1,n-s}(\K)$. 
We know from point (ii) of Paragraph \ref{strategy} that 
$V$ contains a matrix of the form 
$\begin{bmatrix}
[?]_{s \times s} & [?]_{s \times (n-s)} \\
[?]_{(n-s) \times s} & [?]_{(n-s) \times (n-s)} \\
[0]_{1 \times s} & L'
\end{bmatrix}$. Using point (i) above, we deduce that $V$ contains a matrix of the form
$\begin{bmatrix}
[?]_{s \times s} & [0]_{s \times (n-s)} \\
[?]_{(n-s) \times s} & N \\
[0]_{1 \times s} & L'
\end{bmatrix}$
with $N \in \Mat_{n-s}(\K)$. Using again (S3), we find that $N \in U$. Using again point (i) above, 
we conclude that $V$ contains a matrix of the form
$\begin{bmatrix}
[?]_{s \times s} & [0]_{s \times (n-s)} \\
[?]_{(n-s) \times s} & [0]_{(n-s) \times (n-s)} \\
[0]_{1 \times s} & L'
\end{bmatrix}$. 
Therefore:
\begin{itemize}
\item[(iii)] For every $L' \in \Mat_{1,n-s}(\K)$,
the linear subspace $V$ contains a matrix of the form
$\begin{bmatrix}
[?]_{s \times s} & [0]_{s \times (n-s)} \\
[?]_{(n-s) \times s} & [0]_{(n-s) \times (n-s)} \\
[0]_{1 \times s} & L'
\end{bmatrix}$. 
\end{itemize}

We claim that $\calV=i_{n+1,n}(\calW)$ in this reduced situation.
In order to prove it, we need to distinguish between two cases, whether $s=1$ or $s>1$.

\paragraph{The case $s=1$} ${}$ \\
In this case, every matrix of $W$ has zero as first column.
It follows that we may find a matrix $C_1 \in \Mat_{n-1,1}(\K)$ and two linear maps $f : \Mat_{1,n-1}(\K) \rightarrow \Mat_{n-1,1}(\K)$
and $b : \Mat_{1,n-1}(\K) \rightarrow \K$ such that,
for every $\alpha \in \K$ and every pair $(L,L')\in \Mat_{1,n-1}(\K)^2$,
the space $\calV$ contains the matrix
$$\begin{bmatrix}
1+b(L) & L' \\
\alpha \,C_1+f(L) & I_{n-1} \\
\alpha & L
\end{bmatrix}.$$

\begin{claim}
One has $C_1=0$, $b=0$ and $f=0$.
\end{claim}

\begin{proof}
Let $(L,\alpha)\in \Mat_{1,n-1}(\K) \times \K$.
We first show that $M:=\begin{bmatrix}
\alpha \,C_1+f(L) & I_{n-1} \\
\alpha & L
\end{bmatrix}$ is singular only if $\alpha=0$ and $f(L)=0$. \\
Assume indeed that $M$ is singular. Then, by the same line of reasoning as in the proof of Claim \ref{M1in},
the first column of $M$ must equal zero (otherwise $\rk M=n-1$, and we would be able to find $L' \in \Mat_{1,n-1}(\K)$
such that $\begin{bmatrix}
1+b(L) & L'
\end{bmatrix}$ belongs to the row space of $M$,
which would contradict the fact that $\rk \begin{bmatrix}
1+b(L) & L' \\
\alpha \,C_1+f(L) & I_{n-1} \\
\alpha & L
\end{bmatrix}=n$).
It follows that $\alpha=0$ and $\alpha\,C_1+f(L)=0$, and thus $f(L)=0$. \\
Computing the determinant of $M$ with Gaussian elimination,
we deduce that
\begin{equation}\label{imply}
\forall (\alpha,L) \in \K \times \Mat_{1,n-1}(\K), \; \alpha(1-LC_1)=Lf(L) \; \Rightarrow \; \bigl(\alpha=0  \quad \text{and} \quad f(L)=0\bigr).
\end{equation}
As we shall now see, \eqref{imply} is enough to show that $f=0$ and $C_1=0$. \\
For every $L \in \Mat_{1,n-1}(\K)$ such that $LC_1\neq 1$, we may choose an $\alpha \in \K$ such that
$\alpha(1-LC_1)=Lf(L)$, which yields $f(L)=0$.
Notice that $\{L \in \Mat_{1,n-1}(\K) : \; LC_1 \neq 1\}$ spans $\Mat_{1,n-1}(\K)$: this is obvious indeed if $C_1=0$,
and if not, then $\{L \in \Mat_{1,n-1}(\K) : \; LC_1 \neq 1\}$ contains two parallel affine hyperplanes
of $\Mat_{1,n-1}(\K)$ (as $\# \K>2$) and therefore cannot be included in a linear hyperplane of $\Mat_{1,n-1}(\K)$. Since $f$ is linear,
we deduce that $f=0$. \\
If $C_1 \neq 0$, then we may choose $L \in \Mat_{1,n-1}(\K) \setminus \{0\}$ such that $LC_1=1$, in which
case taking $\alpha=1$ in \eqref{imply} yields a contradiction. Therefore $C_1=0$. \\
Finally, if $b \neq 0$, then we may choose $L \in \Mat_{1,n-1}(\K)$ such that $b(L)=-1$,
and, by taking $\alpha=0$, we deduce that $\calV$ contains the rank $n-1$ matrix
$\begin{bmatrix}
0 & [0]_{1 \times (n-1)} \\
[0]_{(n-1) \times 1} & I_{n-1} \\
0 & L
\end{bmatrix}$, a contradiction.
We conclude that $b=0$.
\end{proof}

We have therefore proven that $i_{n+1,n}(\calW) \subset \calV$ hence $\calV=i_{n+1,n}(\calW)$
since these affine spaces have the same dimension.

\paragraph{The case $s>1$} ${}$ \\
In this case, we start by ``cleaning up" the upper left $r \times r$ blocks:

\begin{claim}\label{claim5}
For every $L_2 \in \Mat_{1,n-s}(\K)$, the linear subspace $V$ contains a matrix of the form
$\begin{bmatrix}
[0]_{s \times s} & [0]_{s \times (n-s)} \\
[?]_{(n-s) \times s} & [0]_{(n-s) \times (n-s)} \\
[0]_{1 \times s} & L_2
\end{bmatrix}$.
\end{claim}

\begin{proof}
Let $L_2 \in \Mat_{1,n-s}(\K)$.
We already know that $V$ contains a matrix of the form
$\begin{bmatrix}
A & [0]_{s \times (n-s)} \\
[?]_{(n-s) \times s} & [0]_{(n-s) \times (n-s)} \\
[0]_{1 \times s} & L_2
\end{bmatrix}$ with $A \in \Mat_s(\K)$.
Our goal is to show that $A \in P\Mata_s(\K)$.
Indeed, if we could prove it, then
we would know from (S1) that $V$ contains the matrix
$\begin{bmatrix}
A & [0]_{s \times (n-s)} \\
[0]_{(n+1-s) \times s} & [0]_{(n+1-s)\times (n-s)}
\end{bmatrix}$, and subtracting it with the above matrix would yield the claimed result.

Set
$L_1:=\begin{bmatrix}
1 & 0 & \cdots & 0
\end{bmatrix} \in \Mat_{1,s}(\K)$
and
$L'_1:=\begin{bmatrix}
0 & \cdots & 0 & 1
\end{bmatrix} \in \Mat_{1,s}(\K)$.
By (S2), $V$ contains a matrix of the form
$\begin{bmatrix}
A & [0]_{s \times (n-s)} \\
[?]_{(n-s) \times s} & [0]_{(n-s) \times (n-s)} \\
L'_1 & L_2
\end{bmatrix}$. Let us now write $L_2=\begin{bmatrix}
a_{s+1} & \cdots & a_n
\end{bmatrix}$ and perform the column operations $C_k \leftarrow C_k-a_kC_s$ for $k$ from $s+1$ to $n$
on $\calV$, in order to recover a new affine subspace $\calV'$ of $\Mat_{n+1,n}(\K)$ that is still roughly reduced. Notice that $\calV'$ still has
$\calW$ as core space and that the translation vector space of $\calV'$ contains a matrix of the form
$B=\begin{bmatrix}
A & [0]_{s \times (n-s)} \\
[?]_{(n-s) \times s} & [?]_{(n-s) \times (n-s)} \\
L'_1 & [0]_{1 \times (n-s)}
\end{bmatrix}$ (here, we have used property (S1) in order to simplify the upper-right block). \\
Note that $T(\calV'_1)$ might be different from $T(\calV_1)$. As we shall now prove, this is not the case. \\
Using the line of reasoning that lead to (S2), we
find a new series of row operations $L_k \leftarrow L_k+\mu_k \,L_{n+1}$ for $k$ from $1$ to $s$
so that the new affine subspace $\calV''$ obtained from $\calV'$ satisfies
$T(\calV_1'')=i_{s+1,s}\bigl(I_s+P\Mata_s(\K)\bigr)$. \\
The translation vector space of $T(\calV_1'')$ contains
$\begin{bmatrix}
[0]_{s \times s} \\
L_1
\end{bmatrix}$. Using the above row operations backwards, we deduce that the translation vector space of $T(\calV_1')$
contains
$\begin{bmatrix}
-C_1 & [0]_{s \times (s-1)} \\
1 & [0]_{1 \times (s-1)}
\end{bmatrix} \in \Mat_{s+1,s}(\K)$, where $C_1:=\begin{bmatrix}
\mu_1 \\
\vdots \\
\mu_s
\end{bmatrix}$. Using the initial column operations backwards, we deduce, as $s \geq 2$, that
$T(\calV_1)$ also contains
$\begin{bmatrix}
-C_1 & [0]_{s \times (s-1)} \\
1 & [0]_{1 \times (s-1)}
\end{bmatrix}$. However we had $T(\calV_1)=i_{s+1,s}\bigl(I_s+P\Mata_s(\K)\bigr)$ and hence
$\begin{bmatrix}
-C_1 & [0]_{s \times (s-1)}
\end{bmatrix}$ belongs to $P\, \Mata_s(\K)$. \\
As $P\Mata_s(\K)$ contains no rank $1$ matrix, we deduce that
$\mu_1=\dots=\mu_s=0$, therefore $T(\calV'_1)=i_{s+1,s}(I_s+P\Mata_s(\K))$.
Considering $B$, we deduce that $A \in P\Mata_s(\K)$, QED.
\end{proof}

Combining (S1), (S2) and Claim \ref{claim5}, we find that for every
$(L_1,L_2)\in \Mat_{1,s}(\K) \times \Mat_{1,n-s}(\K)$, the subspace $V$ contains a unique matrix of the form
$$\begin{bmatrix}
[0]_{s \times s} & [0]_{s \times (n-s)} \\
[?]_{(n-s) \times s} & [0]_{(n-s) \times (n-s)} \\
L_1 & L_2
\end{bmatrix}.$$
We deduce that there are linear maps
$\varphi : \Mat_{1,s}(\K) \rightarrow \Mat_{n-s,s}(\K)$ and
$\psi : \Mat_{1,n-s}(\K) \rightarrow \Mat_{n-s,s}(\K)$ such that,
for every $(L_1,L_2)\in \Mat_{1,s}(\K) \times \Mat_{1,n-s}(\K)$,
the linear subspace $V$ contains the matrix
$$\begin{bmatrix}
[0]_{s \times s} & [0]_{s \times (n-s)} \\
\varphi(L_1)+\psi(L_2) & [0]_{(n-s) \times (n-s)} \\
L_1 & L_2
\end{bmatrix}.$$
As in the proof of Proposition \ref{irrprop}, we lose no generality in assuming that
$P=\begin{bmatrix}
\alpha & [0]_{1 \times (s-1)} \\
[?]_{(s-1) \times 1} & Z
\end{bmatrix}$ for some $\alpha \in \K \setminus \{0\}$ and some non-isotropic matrix $Z \in \GL_{s-1}(\K)$.
In this situation, some information on $\varphi$ and $\psi$ may be obtained by using the induction hypothesis:

\begin{claim}\label{claim6}
For every $(L_1,L_2) \in \Mat_{1,s}(\K) \times \Mat_{1,n-s}(\K)$, the last $s-1$ columns of
$\varphi(L_1)$ and of $\psi(L_2)$ are zero.
\end{claim}

\begin{proof}
Denote by $\calV_2$ the affine subspace of $\calV$ consisting of its matrices
with $\begin{bmatrix}
1 & 0 & \cdots & 0
\end{bmatrix}$ as first row.
For every $M \in \calV_2$, write $M=\begin{bmatrix}
1 & [0]_{1 \times (n-1)} \\
[?]_{n \times 1} & Y(M)
\end{bmatrix}$ with $Y(M) \in \Mat_{n,n-1}(\K)$, and note that $\rk Y(M)=n-1$.
Since $W=P\Mata_s(\K) \vee U$ and $\calV$ contains $\begin{bmatrix}
I_n \\
[0]_{1 \times n}
\end{bmatrix}$, we find:
\begin{itemize}
\item[(i)] For every $A \in Z\Mata_{s-1}(\K)\,\vee\, U$,
the affine subspace $Y(\calV_2)$ contains
$\begin{bmatrix}
I_{n-1}+A \\
[0]_{1 \times (n-1)}
\end{bmatrix}$.
\end{itemize}
Using the definition of $\varphi$ and $\psi$, we also find:
\begin{itemize}
\item[(ii)] For every $L \in \Mat_{1,n-1}(\K)$, the translation vector space of $Y(\calV_2)$
contains a matrix with $L$ as last row.
\end{itemize}
Again, we deduce that $\calV_2$ is an affine subspace of $\Mat_{n,n-1}(\K)$
with lower rank $n-1$ and codimension $\binom{n}{2}$. Moreover, $\calV_2$ is roughly-reduced
with $\calW_2=I_{n-1}+\bigl(Z\Mata_{s-1}(\K)\vee U\bigr)$ as core space.
Applying the induction hypothesis to $Y(\calV_2)$ yields a list
$(a_2,\dots,a_n)\in \K^{n-1}$ such that, for the space $\calV'$ obtained from $\calV$ by the series of row operations
$L_k \leftarrow L_k+a_kL_{n+1}$ for $k$ from $2$ to $n$, one has $Y(\calV'_2)=i_{n,n-1}(\calW_2)$. \\
Let us show that $a_{s+1}=\cdots=a_n=0$. \\
In order to do so, we set $L_2:=\begin{bmatrix}
1 & 0 & \cdots & 0
\end{bmatrix} \in \Mat_{1,n-s}(\K)$. Since the translation vector space of $Y(\calV'_2)$
contains $\begin{bmatrix}
[0]_{(n-1) \times (s-1)} & [0]_{(n-1) \times (n-s)} \\
[0]_{1 \times (s-1)} & L_2
\end{bmatrix}$, taking the above row operations backwards shows that the translation vector space of $Y(\calV_2)$ contains
$$\begin{bmatrix}
[0]_{(s-1)\times (s-1)} & -C_1 & [0]_{(s-1) \times (n-s-1)} \\
[0]_{(n-s)\times (s-1)} & -C_2 & [0]_{(n-s) \times (n-s-1)} \\
0 & 1 & [0]_{1 \times (n-s-1)}
\end{bmatrix}$$
where $C_1=\begin{bmatrix}
a_2 \\
\vdots \\
a_s
\end{bmatrix}$
and
$C_2=\begin{bmatrix}
a_{s+1} \\
\vdots \\
a_n
\end{bmatrix}$.
Therefore $\begin{bmatrix}
[0]_{(s-1)\times (s-1)} & -C_1 & [0]_{(s-1) \times (n-s-1)} \\
[0]_{(n-s)\times (s-1)} & -C_2 & [0]_{(n-s) \times (n-s-1)}
\end{bmatrix}$ belongs to the translation vector space of the core space of $Y(\calV_2)$, i.e., to
$P\Mata_s(\K)\vee P_2 \Mata_{n_2}(\K) \vee \cdots \vee P_q \Mata_{n_q}(\K)$.
As $P_2 \Mata_{n_2}(\K)$ contains no rank $1$ matrix, it follows that $C_2=0$.

Let finally $(L_1,L_2) \in \Mat_{1,s}(\K) \times \Mat_{1,n-s}(\K)$.
Since $a_{s+1}=\cdots=a_n=0$, we find that the translation vector space of $\calV'_2$ contains a matrix of the form
$$B=\begin{bmatrix}
[0]_{1 \times s} & [0]_{1 \times (n-s)} \\
[?]_{(s-1) \times s} & [0]_{(s-1) \times (n-s)} \\
\varphi(L_1)+\psi(L_2) & [0]_{(n-s) \times (n-s)} \\
L_1 & L_2
\end{bmatrix}.$$
The matrix obtained from $B$ by deleting its first row and its first column must then belong to
$i_{n,n-1}\bigl(Z\Mata_{s-1}(\K) \vee U\bigr)$, which shows that the last $s-1$ columns of $\varphi(L_1)+\psi(L_2)$
are zero. Since $\varphi(0)=0$ and $\psi(0)=0$, this obviously yields the claimed result.
\end{proof}

\begin{claim}\label{claim7}
One has $\varphi=0$ and $\psi=0$.
\end{claim}

\begin{proof}
Denote by $(e_1,\dots,e_s)$ the canonical basis of $\K^s$.
Let $Q \in \GL_s(\K)$ and assume that
$QPQ^T$ has the same basic shape as $P$, i.e.,
$QPQ^T=\begin{bmatrix}
\alpha' & [0]_{1 \times (s-1)} \\
[?]_{(s-1) \times 1} & Z'
\end{bmatrix}$ for some $\alpha' \in \K \setminus \{0\}$ and some $Z' \in \GL_{s-1}(\K)$.
Then, multiplying $\calV$ on the left by
$Q \oplus I_{n-s+1}$ and on the right by
$Q^{-1} \oplus I_{n-s}$, we find that
the previous situation is essentially unchanged for the new affine subspace $\calV'$,
the only noticeable difference being that
$\calW$ is replaced with $I_n+\Bigl[(QPQ^T)\Mata_s(\K) \vee U\Bigr]$,
whilst $\varphi$ and $\psi$ are replaced respectively with
$L_1 \mapsto \varphi(L_1Q)Q^{-1}$ and $L_2 \mapsto \psi(L_2)Q^{-1}$.

Let $(L_1,L_2) \in \Mat_{1,s}(\K) \times \Mat_{1,n-s}(\K)$. Using Claim \ref{claim6}, we find that $\varphi(L_1)$ and
$\psi(L_2)$ vanish on $Q^{-1}e_i$ for every $i \in \lcro 2,s\rcro$.
However we already know that they completely vanish on the linear hyperplane $\Vect(e_2,\dots,e_s)$.
If $Q$ does not stabilize $\Vect(e_2,\dots,e_s)$, then it follows that $\varphi=0$ and $\psi=0$, as claimed.
It thus remains to show that $Q$ may be chosen as such, which amounts to proving
that it may be chosen so as to have $Q^Te_1$ linearly independent from $e_1$.
However the bilinear form $b: (X,Y) \mapsto X^TPY$ is non-isotropic (i.e., $X \mapsto X^TPX$ is non-isotropic)
hence $e_2$ may be completed to a basis $(e_2,f_2,\dots,f_s)$ of $\K^s$
so that $b(e_2,f_k)=0$ for every $k \in \lcro 2,s\rcro$.
Denoting by $Q_1$ the matrix of coordinates of $(e_2,f_2,\dots,f_s)$ in $(e_1,\dots,e_s)$, it follows
that the matrix of $b$ in $(e_2,f_2,\dots,f_s)$ has $\begin{bmatrix}
b(e_2,e_2) & 0 & \cdots & 0
\end{bmatrix}$ as first row and equals $Q^T_1PQ_1$. However $Q_1e_1=e_2$, therefore $Q:=Q_1^T$ fulfills our needs,
which proves that $\varphi=0$ and $\psi=0$.
\end{proof}

Since $\calV$ is roughly-reduced with $\calW$ as core space, Claim \ref{claim7} entails that
$i_{n+1,n}(\calW) \subset \calV$,
hence $i_{n+1,n}(\calW)=\calV$ as these affine spaces have the same dimension.
This finishes our proof of Proposition \ref{caser=p=n-1}.

\subsection{The existence statement in the general case}\label{generalcasesection}

In order to prove the general case in the existence statement of Theorem \ref{classaffineboundedbelow},
we first need to establish the following lemma:

\begin{lemme}\label{lastlemma}
Let $\calW$ be a maximal affine subspace of non-singular matrices of $\Mat_n(\K)$, with $n \geq 2$.
Let $\calV$ be an affine subspace of $\Mat_{n+1}(\K)$ and assume:
\begin{enumerate}[(i)]
\item That
$\lrk(\calV) \geq n$;
\item That for every $M \in \calW$ and
$(L,C)\in \Mat_{1,n}(\K) \times \Mat_{n,1}(\K)$, the subspace $\calV$ contains the matrix
$$\begin{bmatrix}
M & C \\
L & 0
\end{bmatrix};$$
\item That $\calV$ contains a matrix of the form
$\begin{bmatrix}
[?]_{n \times n} & [?]_{n \times 1} \\
[?]_{1 \times n} & 1
\end{bmatrix}$.
\end{enumerate}
Then the translation vector space $V$ of $\calV$ contains the matrix
$\begin{bmatrix}
[0]_{n \times n} & [0]_{n \times 1} \\
[0]_{1 \times n} & 1
\end{bmatrix}$.
\end{lemme}

\begin{proof}
Denote by $W$ the translation vector space of $\calW$.
The assumptions show that $V$ contains
$\begin{bmatrix}
[0]_{n \times n} & C \\
L & 0
\end{bmatrix}$ for every $(L,C)\in \Mat_{1,n}(\K) \times \Mat_{n,1}(\K)$
and that it contains
$\begin{bmatrix}
A & [0]_{n \times 1} \\
[0]_{1 \times n} & 1
\end{bmatrix}$ for some $A \in \Mat_n(\K)$.
If $A \in W$, then the claimed result follows immediately. \\
Let us perform a \emph{reductio ad absurdum} by assuming that $A \notin W$.
Then $\calW$ is a strict subspace of the affine space $\calW+\Vect(A)$, and
hence some matrix of $\calW+\Vect(A)$ must be singular.
Since no matrix of $\calW$ is singular, $M+A$ is singular for some $M \in \calW$,
and hence $\calV$ contains
$\begin{bmatrix}
B & [0]_{n \times 1} \\
[0]_{1 \times n} & 1
\end{bmatrix}$ for some singular matrix $B \in \Mat_n(\K)$. \\
Note that $\rk B \geq n-1$ since $\lrk(\calV) \geq n$, therefore $\rk B=n-1$.
There exists $(P,Q)\in \GL_n(\K)^2$ such that
$PBQ=\begin{bmatrix}
I_{n-1} & [0]_{(n-1) \times 1} \\
[0]_{1 \times (n-1)} & 0
\end{bmatrix}$.
Setting $\calV':=\widetilde{P}\,\calV\,\widetilde{Q}$, where
$\widetilde{P}:=P \oplus 1 \in \GL_{n+1}(\K)$ and $\widetilde{Q}:=Q \oplus 1  \in \GL_{n+1}(\K)$,
we see that $\calV'$ satisfies properties (i), (ii) ($\calW$ being replaced with $P\,\calW\,Q$) and (iii) and that it contains
$\begin{bmatrix}
I_{n-1} & [0]_{(n-1) \times 1} & [0]_{(n-1) \times 1} \\
[0]_{1 \times (n-1)} & 0 & 0 \\
[0]_{1 \times (n-1)} & 0 & 1
\end{bmatrix}$.
We then define the row matrix $L:=\begin{bmatrix}
1 & 0 & \cdots & 0
\end{bmatrix} \in \Mat_{1,n-1}(\K)$ and note that the matrix
$\begin{bmatrix}
I_{n-1} & [0]_{(n-1) \times 1} & L^T \\
[0]_{1 \times (n-1)} & 0 & 0  \\
L & 0 & 1
\end{bmatrix}$ belongs to $\calV'$: this is a contradiction since,
judging from its first, $n$-th and $(n+1)$-th rows, this matrix has rank lesser than $n$
(this uses the fact that $n \geq 2$).
This finishes our proof.
\end{proof}

With the previous results, we may now conclude our proof of the existence statement in Theorem \ref{classaffineboundedbelow}.
Let $\calV$ be an affine subspace of $\Mat_{n,p}(\K)$ with codimension $\binom{r+1}{2}$
and such that $\lrk(\calV)=r$. We lose no generality in assuming that $\calV$
is roughly reduced. In that case, we define the spaces $V$ and $\calW$ as in Paragraph \ref{strategy}.

\begin{Not}
We denote by $G$ the linear subspace of $V$ consisting of its matrices of the form
$\begin{bmatrix}
[?]_{r \times r} & [?]_{r \times (p-r)} \\
[?]_{(n-r) \times r} & [0]_{(n-r) \times (p-r)}
\end{bmatrix}$.
\end{Not}

\vskip 3mm
Assume first that $n>r$. Define $G'$ as the linear subspace of $G$
consisting of its matrices in which the $p-r$ last columns
and the $n-r-1$ last rows are zero.
We may then write every matrix of $G'$ as
$M=\begin{bmatrix}
K(M) & [0]_{(r+1) \times (p-r)} \\
[0]_{(n-r-1) \times r} & [0]_{(n-r-1) \times (p-r)}
\end{bmatrix}$ with $K(M) \in \Mat_{r+1,r}(\K)$. Using properties
(i) and (ii) in Paragraph \ref{strategy}, one finds that
$\dim K(G')=(r+1)r-\binom{r+1}{2}$. Moreover, with $J:=\begin{bmatrix}
I_r \\
[0]_{1 \times r}
\end{bmatrix}$, the affine subspace
$J+K(G')$ has lower rank $r$, codimension $\binom{r+1}{2}$ in $\Mat_{r+1,r}(\K)$ and is roughly-reduced
with $\calW$ as core space. Using Proposition \ref{caser=p=n-1},
we then find a list $(\lambda_1,\dots,\lambda_r) \in \K^r$
such that the linear subspace deduced from $V$ by the row operations $L_k \leftarrow L_k+\lambda_k\,L_{r+1}$,
for $k$ from $1$ to $r$, contains every rank $1$ matrix of $\Mat_{n,p}(\K)$ with all columns zero starting from the $(p+1)$-th and all rows
zero except the $(r+1)$-th. \\
Release now the assumption $n>r$: for any $i \in \lcro r+1,n\rcro$, swapping the $i$-th and $(r+1)$-th rows
leaves our basic assumptions unchanged, therefore we may find a list
$(\lambda_1^{(i)},\dots,\lambda_r^{(i)}) \in \K^r$ such that the row operations $L_k \leftarrow L_k+\lambda_k^{(i)} L_i$, for $k$ from $1$ to $r$,
transform $\calV$ into an affine subspace whose translation vector space
contains every rank $1$ matrix of $\Mat_{n,p}(\K)$ with all columns zero starting from the $(p+1)$-th and all rows
zero except the $i$-th. \\
Using the row operations $L_k \leftarrow L_k+\underset{i=r+1}{\overset{n}{\sum}} \lambda_k^{(i)} L_i$,
for $k$ from $1$ to $r$, we thus reduce the situation to the one where:
\begin{center}
(T1) : \quad $G$ contains the matrix
$\begin{bmatrix}
[0]_{r \times r} & [0]_{r \times (p-r)} \\
B & [0]_{(n-r) \times (p-r)}
\end{bmatrix}$ for every $B \in \Mat_{n-r,r}(\K)$.
\end{center}
Note that the above row operations do not change the fact that $\calV$ is roughly-reduced with $\calW$ as core space.

\vskip 3mm
Applying the same technique to $\calV^T$, we find lists $(\mu_1^{(j)},\dots,\mu_r^{(j)}) \in \K^r$, for $j \in \lcro r+1,p\rcro$,
such that the column operations $C_k \leftarrow C_k+\underset{j=r+1}{\overset{p}{\sum}}\mu_k^{(j)} C_j$, for $k$ from $1$ to $r$,
further transform $\calV$
into an affine subspace whose associated linear subspace $G$ contains the matrix
$\begin{bmatrix}
[0]_{r \times r} & C \\
[0]_{(n-r)\times r} & [0]_{(n-r) \times (p-r)}
\end{bmatrix}$ for every $C \in \Mat_{r,p-r}(\K)$. \\
Notice that those column operations preserve the $n-r$ last rows of the matrices of $G$,
and hence property (T1) is still satisfied after performing them
(and $\calV$ remains roughly-reduced with $\calW$ as core space).
We have reduced the situation to the one where :
\begin{center}
(T2) : \quad $G$ contains
$\begin{bmatrix}
[0]_{r \times r} & C \\
B & [0]_{(n-r) \times (p-r)}
\end{bmatrix}$ for every $B \in \Mat_{n-r,r}(\K)$ and every $C \in \Mat_{r,p-r}(\K)$.
\end{center}
\noindent In this reduced situation, we claim that $\calV=i_{n,p}(\calW)$. This is immediate indeed if $n=r$ or $p=r$.
Assume now that $n>r$ and $p>r$.
Using the equality of dimensions between $i_{n,p}(\calW)$ and $\calV$ together with
property (T2), it obviously suffices to show that $V$
contains
$\begin{bmatrix}
[0]_r & [0]_{r \times (p-r)} \\
[0]_{(n-r) \times r} & D
\end{bmatrix}$ for every $D \in \Mat_{n-r,p-r}(\K)$. \\
For $(i,j)\in \lcro 1,n-r\rcro \times \lcro 1,p-r\rcro$, denote by $E_{i,j}\in \Mat_{n-r,p-r}(\K)$
the elementary matrix with entry $1$ at the spot $(i,j)$ and zero elsewhere. \\
Denote by $\calE$ the affine space of all matrices of $\calV$ with the form
$$M=\begin{bmatrix}
\Delta(M) & [0]_{(r+1) \times (p-r-1)} \\
[0]_{(n-r-1) \times (r+1)} & [0]_{(n-r-1)\times (p-r-1)}
\end{bmatrix} \quad \text{for some $\Delta(M) \in \Mat_{r+1}(\K)$.}$$
Using property (T2), we find that the affine subspace $\Delta(\calE)$ of $\Mat_{r+1}(\K)$
satisfies the assumptions of Lemma \ref{lastlemma}.
It follows that $V$ contains the matrix
$\begin{bmatrix}
[0]_{r \times r} & [0]_{r \times (p-r)} \\
[0]_{(n-r) \times r} & E_{1,1}
\end{bmatrix}$. However, our assumptions on $\calV$ are unchanged by using
a row swap $L_i \leftrightarrow L_j$ for an arbitrary $(i,j)\in \lcro r+1,n\rcro^2$
and/or a column swap $C_i \leftrightarrow C_j$ for an arbitrary $(i,j)\in \lcro r+1,p\rcro^2$.
We deduce that the linear subspace $V$ contains
$\begin{bmatrix}
[0]_{r \times r} & [0]_{r \times (p-r)} \\
[0]_{(n-r) \times r} & E_{i,j}
\end{bmatrix}$ for every $(i,j)\in \lcro 1,n-r\rcro \times \lcro 1,p-r\rcro$,
and it thus contains
$\begin{bmatrix}
[0]_{r \times r} & [0]_{r \times (p-r)} \\
[0]_{(n-r) \times r} & D
\end{bmatrix}$ for every $D \in \Mat_{n-r,p-r}(\K)$. This completes the proof of Theorem \ref{classaffineboundedbelow}.

\end{document}